\documentclass[12pt]{amsart}
\usepackage{amscd,amssymb}
\usepackage[arrow,matrix]{xy}

\newtheorem{theorem}[subsection]{Theorem}
\newtheorem{lemma}[subsection]{Lemma}
\newtheorem{prop}[subsection]{Proposition}
\newtheorem{corollary}[subsection]{Corollary}

\newtheorem{thm}{Theorem}

\theoremstyle{definition}
\newtheorem{remark}[subsection]{Remark}

\newtheorem{example}[subsection]{Example}
\newtheorem*{ack}{Acknowledgment}

\newcommand{\F}{{\mathcal F}}

\newcommand{\V}{{\mathcal V}}

\newcommand{\sV}{{\sf V}}
\newcommand{\sE}{{\sf E}}

\newcommand{\T}{\mathbb{T}}
\newcommand{\Z}{\mathbb{Z}}

\newcommand{\RR}{\mathbb{R}}
\newcommand{\C}{\mathbb{C}}

\newcommand{\PP}{\mathbb{P}}

\newcommand{\D}{\mathbb{D}}

\newcommand{\tX}{\widetilde{X}}
\newcommand{\tE}{\widetilde{E}}
\newcommand{\tH}{\widetilde{H}}
\newcommand{\tM}{\widetilde{M}}
\newcommand{\hX}{\widehat{X}}
\newcommand{\hh}{\hat{h}}

\DeclareMathOperator{\Hom}{Hom}
\DeclareMathOperator{\rank}{rank}
\DeclareMathOperator{\im}{im}

\DeclareMathOperator{\tor}{Tor}

\DeclareMathOperator{\genus}{genus}

\DeclareMathOperator{\cd}{cd}
\DeclareMathOperator{\ann}{Ann}

\newcommand{\abs}[1]{\left| #1 \right|}

\newenvironment{romenum}
{

\begin{enumerate}}{\end{enumerate}}

\providecommand{\bigtimes}{\mathop{%
\mathchoice{\raisebox{-2pt}{\huge$\times$}}{\mbox{\Large$\times$}}%
{\raisebox{0pt}{\Large$\times$}}{\times}}\displaylimits}%

\begin{document}

\title[Non-finiteness properties of projective groups]%
{Non-finiteness properties of fundamental groups 
of smooth projective varieties}

\author[A. Dimca]{Alexandru Dimca}
\address{  Laboratoire J.A.~Dieudonn\'{e}, UMR du CNRS 6621, 
                 Universit\'{e} de Nice-Sophia-Antipolis,
                 Parc Valrose,
                 06108 Nice Cedex 02,
                 France}
\email
{dimca@math.unice.fr}

\author[\c{S}. Papadima ]{\c{S}tefan Papadima$^{1}$}
\address{Inst. of Math. ``Simion Stoilow",
P.O. Box 1-764,
RO-014700 Bucharest, Romania}
\email
{Stefan.Papadima@imar.ro}

\author[A.~Suciu]{Alexander~I.~Suciu$^{2}$}
\address{Department of Mathematics,
Northeastern University,
Boston, MA 02115, USA}
\email{a.suciu@neu.edu}

\thanks{$^1$Partially supported by the CEEX Programme 
of the Romanian Ministry of Education and Research,
contract 2-CEx 06-11-20/2006}
\thanks{$^2$Partially supported by NSF grant DMS-0311142}




\begin{abstract}
For each integer $n\ge 2$, we construct an irreducible, 
smooth, complex projective variety $M$ of dimension $n$,
whose fundamental group has infinitely generated homology 
in degree $n+1$ and whose universal cover is a Stein manifold, 
homotopy equivalent to an infinite bouquet of $n$-dimensional 
spheres. This non-finiteness phenomenon is also reflected 
in the fact that the homotopy group $\pi_n(M)$, viewed as a 
module over $\Z \pi_1(M)$, is free of infinite rank.  As a result, 
we give a negative answer to a question of Koll\'{a}r on the 
existence of quasi-projective classifying spaces (up to 
commensurability) for the fundamental groups of smooth 
projective varieties. To obtain our examples, we develop 
a complex analog of a method in geometric group theory 
due to Bestvina and Brady.
\end{abstract}

\maketitle

\section{Introduction and statement of results}
\label{sec:intro}

\subsection{}
\label{ss11}

Let $M$ be an irreducible, smooth complex projective variety, 
with fundamental group $G=\pi_1(M)$. Groups $G$ arising in this 
fashion are called {\em projective groups}. Except for the obvious 
fact that projective groups are finitely presentable, very little 
is known about their finiteness properties.  The aim of this note 
is to answer several questions in this direction.

The classical finiteness conditions in group theory that we 
have in mind are: 

\medskip\noindent
(i) Wall's property $\F_n$ ($n\le \infty$), requiring the 
existence of a classifying space $K(G,1)$ with finite $n$-skeleton. 
Note that $\F_1$ is equivalent to finite generation, whereas 
$\F_2$ is equivalent to finite presentability of $G$.

\medskip\noindent
(ii) Property $FP_n$ ($n\le \infty$), 
requiring the existence of a projective $\Z G$-resolution 
of the trivial $G$-module $\Z$, which is finitely generated 
in all dimensions $\le n$. Note that the $FP_n$ condition 
implies the finite generation of the homology groups 
$H_i(G, \Z)$, for all $i\le n$.  Clearly, if $G$ is of type 
$\F_n$, then $G$ is of type $FP_n$. It follows from \cite{W} 
that the converse also holds, provided $G$ is finitely presentable, 
and $n<\infty$.  

\medskip\noindent
(iii) Finiteness conditions related to the cohomological dimension 
of $G$, defined as $\cd(G)=\sup \,\{\, i \mid H^i(G, A)\ne 0\}$, 
where $A$ ranges over all $\Z G$-modules. 

\medskip
We refer to the works of Bieri \cite{Bi}, Brown \cite{B}, and  
Serre \cite{Se}  for detailed information on finiteness properties 
of groups.

\subsection{}
\label{ss12}

The first question we consider here was formulated
by J.~Koll\'{a}r in \cite[\S0.3.1]{K}: {\em Is a projective group 
$G$ commensurable (up to finite kernels) with another group $G'$, 
admitting a $K(G', 1)$ which is a quasi-projective variety}? 

Note that necessarily $G'$ must have a finite $K(G', 1)$, 
since every quasi-projective variety has the homotopy type
of a finite CW-complex, see  e.g.~\cite[p.~27]{D}.  For a 
discussion of various commensurability notions, we refer 
to \S\ref{subs:commens}.

The second problem we examine here is related to the 
Shafarevich conjecture \cite{S}, as reformulated in geometric 
finiteness terms by Koll\'{a}r, in \cite[0.3.1.1--0.3.1.2]{K}:  
{\em What other kind of finiteness properties are imposed 
on the group $G=\pi_1(M)$ by the Stein property of the 
universal cover, $\widetilde{M}$}? 

Recall that a {\em Stein manifold}\/ is a complex manifold 
which can be biholomorphically embedded as a closed 
subspace of some affine space $\C^r$. A classical result 
of Andreotti and Frankel \cite{AF} asserts that a Stein 
manifold of (complex) dimension $n$ has the homotopy 
type of a CW-complex of dimension at most $n$. 

\subsection{}
\label{ss13}

The first example of a finitely presented group with infinitely 
generated third homology group is due to J.~Stallings \cite{Sta}. 
A systematic way of constructing groups $N$ of type $\F_n$, 
but not of type $FP_{n+1}$, was found by M.~Bestvina and 
N.~Brady \cite{BB}.  These authors start with a finite graph 
$\Gamma=(\sV,\sE)$, and consider the associated right-angled 
Artin group $G_{\Gamma}$, with a generator $v$ for each 
vertex in $\sV$, and with a relation $uv=vu$ for each edge 
in $\sE$. The {\em Bestvina--Brady group}\/ $N_{\Gamma}$ 
is then defined as the kernel of the homomorphism 
$\nu \colon G_{\Gamma}\to \Z$, which sends each 
generator $v$ to $1$.

The group $G_{\Gamma}$ admits as classifying space a 
subcomplex $K_{\Gamma}$ of the torus of dimension 
$\abs{\sV}$, with cells in one-to-one correspondence 
with the simplices of the flag complex $\Delta_{\Gamma}$. 
Bestvina and Brady had the remarkable idea of exploiting 
the natural affine cell structure of the universal 
cover, $\widetilde{K_{\Gamma}}$, and to do a geometric 
and combinatorial version of real Morse theory on 
$\widetilde{K_{\Gamma}}$. In this way, they were able 
to establish a spectacular connection between the finiteness 
properties of the group $N_{\Gamma}$, and the homotopical 
properties of the simplicial complex $\Delta_{\Gamma}$.

It was noticed in \cite{PS} that the Stallings group may be 
realized as the fundamental group of the complement of a 
complex line arrangement in $\PP^2$. In \cite{DPS06}, we 
identified a large class of Bestvina--Brady groups which 
are {\em quasi-projective}, yet are not commensurable 
to any group admitting a classifying space which is a 
quasi-projective variety.  Starting from a group 
$G_{\Gamma}$ which is a product of $r\ge 3$ free 
groups on at least two generators, we showed that 
$N_{\Gamma}=\pi_1(H)$, where $H$ is the generic 
fiber of an explicit polynomial map, 
$h\colon X\to \C^*$, and $X=\C^r \setminus \{ h=0\}$. 
Thus, $N_{\Gamma}$ is the fundamental group 
of an irreducible, smooth complex affine variety 
of dimension $r-1$. On the other hand, 
$H_r(N_{\Gamma};\Z)$ is not finitely generated, 
and so $N_{\Gamma}$ is not of type $FP_r$. 

\subsection{}
\label{ss14}

In this paper, we develop a complex analog of the Bestvina--Brady 
method, well adapted to construct projective groups with 
controlled finiteness properties. This leads to answers 
to the two questions mentioned in \S\ref{ss12}, as follows.

\begin{thm}
\label{thm:exko}
For each $r\ge 3$, there is an $(r-1)$-dimensional, smooth, 
irreducible, complex projective variety $H$, with fundamental 
group $G$, such that:
\begin{enumerate}
\item \label{ko1}
The  homotopy groups $\pi_i(H)$ vanish for $2\le i\le r-2$, 
while $\pi_{r-1}(H)\ne 0$.
\item \label{ko0}
The universal cover $\widetilde{H}$ is a Stein manifold.
\item \label{ko2}
The group $G$ is of type $\F_{r-1}$, but not of type $FP_r$.
\item \label{ko3}
The group $G$ is not commensurable (up to finite kernels) 
to any group having a classifying space of finite type.
\end{enumerate}
\end{thm}

Our method actually yields stronger results (see 
Corollary \ref{cor:hipi}): the first non-vanishing 
higher homotopy group $\pi_{r-1}(H)$ from \eqref{ko1} is a free 
$\Z \pi_1(H)$-module, with geometrically computable system 
of free generators, and the universal cover $\widetilde{H}$ 
in \eqref{ko0} has the homotopy type of a wedge of 
$(r-1)$-spheres.  Similar results were obtained in \cite{DP}, 
for {\em open}, smooth algebraic varieties $H$.

Theorem \ref{thm:exko} gives a negative answer to Koll\'{a}r's 
question: by Part \eqref{ko3}, the group $G$ is not 
commensurable (up to finite kernels) to any group $G'$ 
admitting a $K(G', 1)$ which is a quasi-projective variety. 

As for the second question, it is easy to show that the 
Stein property of the universal cover of a smooth projective 
variety $M$ forces the cohomological dimension of the 
projective group $G=\pi_1(M)$ to be  larger than the 
complex dimension of $M$; see Proposition \ref{prop:cdim}.
On the other hand, there is no implication of the Stein 
condition at the level of $FP_*$ finiteness properties of $G$.  
To see this, compare smooth projective curves $C$ of positive 
genus, for which $\widetilde{C}$ is Stein and $C$ is aspherical, 
to the varieties $H$ in Theorem \ref{thm:exko}, for which 
$\widetilde{H}$ is Stein, yet $\pi_1(H)$ is not of type $FP_r$, 
for $r=\dim H +1$.

Finally, let us note that Theorem \ref{thm:exko} also sheds 
light on the following question of Johnson and Rees \cite{JR}: 
are fundamental groups of compact K\"{a}hler manifolds 
Poincar\'{e} duality groups of even cohomological dimension? 
In \cite{To}, Toledo answered this question, by producing examples 
of smooth projective varieties $M$ with $\pi_1(M)$ of {\em odd}\/ 
cohomological dimension. Our results (see Theorem  
\ref{thm:exo}\eqref{e2}) show that  fundamental 
groups of smooth projective varieties need not be Poincar\'{e} 
duality groups of {\em any}\/ cohomological dimension:   their 
Betti numbers need not be finite. 

\subsection{}
\label{ss15}

We start by establishing a general relationship between the 
finiteness properties of proper normal subgroups $N$ of a 
finitely generated group $G$, with $G/N$ torsion-free abelian, 
and the structure of certain subsets of the complex algebraic 
torus $\T_G =\Hom (G, \C^*)$.  

By definition, the {\em characteristic varieties}\/ of $G$ are the 
jumping loci for the homology of $G$ with coefficients in rank 
one complex local systems:
\begin{equation}
\label{eq:defv}
\V^s_t(G)= \{ \rho\in \T_G  \mid  \dim_{\C} H_s(G, \C_{\rho})\ge t \}\, .
\end{equation} 
If $G$ is of type $FP_n$, it is readily seen that $\V^s_t(G)$ is an 
algebraic subvariety of $\T_G$, for $s<n$. If $G$ is finitely 
presented, the varieties $\V^1_t(G)$ may be computed directly 
from a presentation of $G$, using the Fox free differential calculus, 
see e.g.~\cite{Hi}.

The importance of these varieties emerged from 
work of S.P.~Novikov \cite{N} on Morse theory for closed 
$1$-forms on manifolds. As shown by Arapura \cite{A}, 
the characteristic varieties $\V^1_1(G)$ provide powerful 
obstructions for deciding the realizability of $G$ as the 
fundamental group of a smooth quasi-projective variety.  
We refer to \cite{DPS05} for various refinements in the 
case of $1$-formal groups, in particular, projective groups.

\begin{thm}
\label{thm:vfin}
Let $G$ be a finitely generated group. 
Suppose $\nu\colon G\to A$ is a non-trivial homomorphism 
to a torsion-free abelian group $A$, and set $N= \ker(\nu)$.  
If $\V^r_1(G)=\T_G$ for some integer $r \ge 1$, then:
\begin{enumerate}
\item \label{b1}
$\dim_{\C} H_{\le r}(N, \C)=\infty$.
\item \label{b2}
$N$ is not commensurable (up to finite kernels) to any 
group of type $FP_r$. 
\end{enumerate}
\end{thm}

The proof is given in Section \ref{sec:charfin}.  
The theorem applies to groups of the form 
$G=\bigtimes_{j=1}^{r} \pi_1(C_j)$, with each $C_j$ 
a smooth complex curve of negative Euler characteristic.

\subsection{}
\label{ss16}

As noted by Deligne \cite{De}, every finitely presentable 
group can be realized as the fundamental group of 
an algebraic variety $X$ (which can be chosen 
as the union of an arrangement of affine subspaces 
in some $\C^n$).  Insisting that $X$ be a {\em smooth}\/ 
variety, though, puts some stringent conditions on what 
groups can occur;  see the monograph \cite{ABC} as 
a reference, and \cite{DPS05} for some recent developments.

We set out here to construct new, interesting examples 
of projective groups.  With Theorem \ref{thm:vfin} in mind, 
we adopt the viewpoint of realizing these new groups as 
subgroups of known groups, rather than extensions of 
known groups.  
 
A convenient setup is provided by {\em irrational pencils}, 
that is, holomorphic maps $h\colon X\to E$ between compact, 
connected, complex analytic manifolds, with target a curve 
$E$ with $\chi (E)\le 0$, and with connected generic fiber.  
Let $p\colon \widetilde{E}\to E$ be the universal cover, and 
denote by $\hh\colon \hX\to \widetilde{E}$ the pull-back of 
$h$ via $p$. Clearly, the maps $h$ and $\hh$ have 
the same fibers; let $H$ be the common {\em smooth}\/ fiber.
Using complex Morse theory on $\hX$, we obtain the following. 

\begin{thm}
\label{thm:cplxmorse}
Let $h\colon X\to E$ be an irrational pencil. 
Suppose $h$ has only isolated singularities. Then:

\begin{enumerate}
\item \label{c1} 
$\pi_i(\hX, H)=0$, for all $i<\dim X$. 

\item \label{c2}
If, moreover, $\dim X\ge 3$, the induced homomorphism 
$h_{\sharp} \colon \pi_1(X)\to \pi_1(E)$ is surjective, 
with kernel isomorphic to $\pi_1(H)$. 
\end{enumerate}
\end{thm}

The proof is given in Section \ref{sec:bbcplx}.  For more on 
complex Morse theory, see Looijenga's book \cite{L}, 
especially Section 5B.  

A concrete family of examples is constructed   
in Section \ref{sect:elliptic}. Starting with an elliptic 
curve $E$, we take $2$-fold branched covers 
$f_j\colon C_j \to E$ ($1\le j \le r$), so that each 
curve $C_j$ has genus at least $2$.  Setting 
$X=\prod_{j=1}^r C_j$, 
we see that $X$ is a smooth, projective variety, 
whose universal cover is a contractible, Stein manifold.   
Moreover, $\V^r_1(\pi_1(X))=\T_{\pi_1(X)}$.  
Using the group law on $E$, we can define a 
map $h\colon X\to E$ by $h=\sum_{j=1}^r f_j$, for $r\ge 2$.
Under certain assumptions on the branched covers 
$f_j$, we show that the smooth fiber of $h$ is 
connected, and $h$ has only isolated  singularities. 

Invoking now Theorems \ref{thm:vfin} and \ref{thm:cplxmorse} 
completes the proof of Theorem \ref{thm:exko}.  
Details and further discussion are given in 
Section \ref{sec:exotic}.

\section{Characteristic varieties and finiteness properties of subgroups}
\label{sec:charfin}

This section is devoted to the proof of Theorem \ref{thm:vfin}.

\subsection{} 
\label{ss21}
Fix a positive integer $m$.  Let $\T_m=\Hom(\Z^m,\C^*)$ be 
the character torus of $\Z^m$, and let $\Lambda= \C \Z^m$ 
be its coordinate ring. Note that $\Lambda$ is isomorphic to 
the ring of Laurent polynomials in $m$ variables.   In particular, 
$\Lambda$ is a noetherian ring, of dimension $m$.

\begin{lemma}
\label{lem:torfdim}
Let $A$ be a $\Lambda$-module which is finite-dimensional as 
a $\C$-vector space.  Then,  for each $j\ge 0$, the set 
\begin{equation*}
\label{eq:zopen}
A_j:=\{ \rho\in \T_m  \mid  \tor^{\Lambda}_j (\C_{\rho}, A)=0 \}
\end{equation*}
is a Zariski open, non-empty subset of the algebraic torus $\T_m$.
\end{lemma}

\begin{proof}
Pick a free $\Lambda$-resolution 
$F_{\bullet}\xrightarrow{\varepsilon} A$, 
with $F_j=\Lambda^{c_j}$, and view the differentials 
$d_j\colon \Lambda^{c_j} \to \Lambda^{c_{j-1}}$ as 
matrices with entries in $\Lambda$.  For a character 
$\rho\in \T_m$, let 
$d_{j}(\rho)\colon \C^{c_j} \to \C^{c_{j-1}}$ be 
the evaluation of $d_j$ at $\rho$. Clearly, 
$\rho\in A_j$ if and only if 
\[
\rank d_{j+1}(\rho) +\rank d_{j}(\rho)\ge c_j,
\] 
a Zariski open condition on $\rho$.  Assuming 
$A_j$ to be empty, we derive a contradiction, as follows.

Let $f\in \ann_{\Lambda}(A)$. Denote by 
$\mu_f$ the homothety induced by $f$ on $\Lambda$-modules. 
Since $\mu_f=0$ on $A$, $\mu_f$ must induce the zero map on 
$\tor^{\Lambda}_j (\C_{\rho}, A)$, for any $\rho\in \T_m$.
In turn, this map may be computed by using $\mu_f$ on 
$F_{\bullet}$. It follows that the homothety $\mu_{f(\rho)}$ 
on $\C^{c_j}$ induces the zero map on $\tor^{\Lambda}_j (\C_{\rho}, A)$, 
which is a non-zero $\C$-vector space, by our assumption 
on $A_j$. Hence, $f(\rho)=0$ for any $\rho \in \T_m$, 
i.e., $f=0$ in $\Lambda$. This shows that $\ann_{\Lambda}(A)=0$.

Now recall that $\dim_{\C} A< \infty$, which implies that the 
$\Lambda$-module $A$ is both noetherian and artinian, therefore 
of finite length. So, $\dim (\Lambda\slash\!\ann_{\Lambda}(A) )=0$, 
by standard commutative algebra. But $\dim (\Lambda)=m>0$, 
a contradiction.
\end{proof}

\subsection{} 
\label{ss23}

Let $G$ be a finitely generated group, and suppose 
$\nu \colon G \to \Z^m$ is an epimorphism. 
Writing  $N=\ker(\nu)$, we have an exact sequence
\begin{equation} 
\label{eq:nexsq}
\xymatrix{
1\ar[r] & N\ar[r] & G \ar^(.48){\nu}[r] & \Z^m \ar[r]& 0}.
\end{equation}
Denote by $\nu^*\colon \T_m=\Hom(\Z^m,\C^*) \to 
\T_G=\Hom(G,\C^*)$ 
the induced map between character tori.

\begin{theorem}
\label{thm:zvfin}
Assume that $\dim_{\C} H_{\le r}(N, \C)<\infty$. 
Then there is a Zariski open, non-empty subset $U\subset \T_m$ 
such that $H_{\le r}(G, \C_{\nu^* \rho})=0$, for any $\rho \in U$.
\end{theorem}

\begin{proof}
By Shapiro's Lemma, $H_*(N, \C)= H_*(G, \Lambda)$.  
Let us examine  the spectral sequence 
associated to the base change $\rho\colon \Lambda \to \C$, 
for a fixed character $\rho\in \T_m$:
\[
E^2_{st}= \tor^{\Lambda}_s (\C_{\rho}, H_t(G, \Lambda)) 
\Rightarrow H_{s+t}(G, \C_{\nu^* \rho})\, ,
\]
see \cite[Theorem XII.12.1]{ML}.
A finite number of applications of Lemma \ref{lem:torfdim} guarantees 
the existence of a Zariski open, non-empty subset $U\subset \T_m$, 
such that $E^2_{st}$ vanishes, provided $s,t\le r$, and $\rho \in U$. 
The conclusion follows.
\end{proof}

\begin{remark}
\label{rem:nov}
If $G$ admits a finite $K(G, 1)$, Theorem \ref{thm:zvfin} also 
follows from Novikov-Morse theory; see \cite{F}, Proposition~1.30 
and Theorem~1.50. See also \cite{DF}, Theorem~1 for a related 
result, under the same finiteness assumption on $G$.
\end{remark}

\subsection{}
\label{subs:commens}

Two groups, $G$ and $G'$, are said to be {\em commensurable}\/ 
if there is a group $\pi$ and a diagram 
\begin{equation*}
\label{eq:commens}
\xymatrixrowsep{12pt}
\xymatrixcolsep{10pt}
\xymatrix{ G & & G' \\ & \pi \ar[ur] \ar[ul] & }
\end{equation*}
with arrows injective and of cofinite image.  The 
two groups are said to be {\em commensurable up to finite kernels}\/ 
if there is a zig-zag of such diagrams, connecting $G$ to $G'$, 
with arrows of finite kernel and cofinite image. Commensurability 
implies commensurability up to finite kernels, but the converse is 
not true in general.  Nevertheless, the two notions coincide if 
one of the two groups is residually finite. For details on all this, 
see the book by de la Harpe \cite[\S IV.B.27--28]{dlH}.

\begin{prop}
\label{prop:commens}
Suppose $G$ and $G'$ are commensurable up to finite kernels. 
Then $G$ is of type $FP_n$ if and only $G'$ is of type $FP_n$. 
\end{prop}

This is an immediate consequence of the following two results 
of Bieri. 

\begin{lemma}[\cite{Bi}, Proposition~2.5]
\label{lem:bieri1}
Let $\pi$ be a finite-index subgroup of $G$. 
Then $G$ is of type $FP_n$ if and only if $\pi$ is.
\end{lemma} 

\begin{lemma}[\cite{Bi}, Proposition~2.7]
\label{lem:bieri2}
Let $1\to N\to G \to Q \to 1$ be an exact sequence of groups, 
and assume $N$ is of type $FP_{\infty}$.  Then $G$ is of type 
$FP_n$ if and only if $Q$ is.
\end{lemma} 

\subsection{} 
\label{subs:proof thm B}

We are now in position to finish the proof of Theorem \ref{thm:vfin}. 
Recall we are given a finitely generated group $G$, and 
a non-trivial homomorphism $\nu\colon G\to A$ to a 
torsion-free abelian group $A$.  Write $N=\ker(\nu)$. 
Note that $\im(\nu)\cong \Z^m$, for some $m>0$.  
Without loss of generality, we may assume $\nu$ is surjective, 
so that we have the exact sequence \eqref{eq:nexsq}.

Part \eqref{b1}.  By assumption, there is an integer $r>0$ 
such that $\V^r_1(G)=\T_G$; that is to say, $H_r(G,\C_{\rho})\ne 0$, 
for all $\rho\in \T_G$.  By Theorem~\ref{thm:zvfin}, 
$\dim_{\C} H_{\le r}(N, \C)=\infty$. 

Part \eqref{b2}. By Part \eqref{b1}, the group $N$ is not 
of type $FP_r$. The conclusion follows from 
Proposition \ref{prop:commens}.   \hfill\qed

\subsection{} 
\label{subs:prodcurves}

We conclude this section with some simple examples of 
groups $G$ to which Theorem \ref{thm:vfin} applies.

\begin{lemma}
\label{lem:vprod}
Let $G=\bigtimes_{i=1}^r G_i$ be a product of finitely 
generated groups. If $\V^1_1(G_i)=\T_{G_i}$, for all $i$, 
then $\V^r_1(G)=\T_G$.
\end{lemma}

\begin{proof}
For a character $\rho\in \T_G$, denote by $\rho_i \in \T_{G_i}$ 
the restriction of $\rho$ to $G_i$. By the K\"{u}nneth formula, 
$H_r(G, \C_{\rho})\supset \bigotimes_{i=1}^r H_1(G_i, \C_{\rho_i})$, 
and the tensor product is non-zero, by hypothesis.
\end{proof}

\begin{example}
\label{ex:curves}
Let $G$ be the fundamental group of a smooth (not necessarily 
compact) complex curve $C$, with $\chi(C)<0$.  Then 
$\V^1_1(G)=\T_G$. Indeed, for any $\rho\in \T_G$, 
the Euler characteristic $\chi(C,\C_{\rho}):=
\dim H_0(C,\C_{\rho})- \dim H_1(C,\C_{\rho})+
\dim H_2(C,\C_{\rho})$ equals $\chi(C)$, and 
the claim follows.
\end{example}

Using Lemma \ref{lem:vprod}, Example \ref{ex:curves}, and 
Theorem \ref{thm:vfin}, we obtain the following. 

\begin{corollary}
\label{cor:product of curves}
Let $C_1,\dots,C_r$ be smooth, complex 
curves with $\chi(C_j)<0$, and let 
$G=\pi_1(\prod_{j=1}^{r} C_j)$. 
Then:
\begin{enumerate}
\item  \label{cc1} 
$\V^r_1(G)=\T_G$.   
\item \label{cc2}
If $N$ is a normal subgroup of $G$, with $G/N\cong \Z^m$, 
for some $m>0$, then $N$ is not commensurable 
(up to finite kernels) to any group of type $FP_r$.
\end{enumerate}
\end{corollary}

In this context, we should note that 
the $FP_n$ finiteness and non-finiteness properties 
of subgroups of finite products of surface groups were 
analyzed by Bridson, Howie, Miller, and Short \cite{BHMS}, 
using different methods. The fact that the subgroup $N$ from
Corollary \ref{cor:product of curves}\eqref{cc2} above cannot be
of type $FP_r$ may be deduced from the results in \cite{BHMS}; 
our Theorem \ref{thm:vfin} improves this to $\dim_{\C} H_{\le r}(N, \C)=\infty$.

\begin{remark}
\label{rem:resfin}
As is well-known, fundamental groups of smooth complex 
curves are residually finite.  Hence, the product groups 
$G$ (and thus, their subgroups $N$) from above 
are also residually finite.  Consequently, the 
two notions of commensurability are equivalent for 
such groups. Note however that there do exist projective 
groups which are not residually finite, see \cite{ABC}. 
\end{remark}

\section{A complex analog of Bestvina--Brady theory}
\label{sec:bbcplx}

In this section, we prove Theorem \ref{thm:cplxmorse} from 
the Introduction.  

\subsection{}
\label{subs:smooth fibers}
Let $X$ and $E$ be compact, connected, complex analytic 
manifolds. Assume $r:=\dim X>1$ and $\dim E=1$.  
Let $h\colon X\to E$ be a holomorphic map, and   
denote by $C(h)$ the set of critical points of $h$. 
Write $E^*=E \setminus h(C(h))$ and $X^*=h^{-1}(E^*)$. 

Since $h$ is a proper map, the restriction 
$h^*\colon X^* \to E^*$ is a topologically locally trivial fibration. 
The fibers $H_t=h^{-1}(t)$, with $t\in E^*$, are called 
the {\em smooth fibers}\/ of $h$. Clearly, such fibers are 
homeomorphic to each other. The map $h$ is called an 
{\em irrational pencil}\/ if $E$ is a curve of positive genus, 
and the smooth fiber of $h$, denoted $H$, is connected.  

So let $h\colon X\to E$ be an irrational pencil, and  
consider the exact homotopy sequence of 
the fibration $H \hookrightarrow X^* \xrightarrow{h^*} E^*$. 
Since $H$ is connected, $h^*_{\sharp}$ is an epimorphism.
Clearly, the inclusion $\iota \colon E^* \to E$ induces an  
epimorphism $\iota_{\sharp} \colon \pi_1( E^*) \to \pi_1(E)$. 
It follows that $h_{\sharp}$ is an epimorphism 
as well. Hence, $h_{\sharp}$ induces an isomorphism
\begin{equation}
\label{eq:hsharp}
\xymatrix{\pi_1(X)/\ker(h_{\sharp}) \ar^(.6){\cong}[r] & \pi_1(E)}.
\end{equation}

\subsection{}
\label{subs:univ cover}
Now let $p\colon \tE \to E$ be the universal covering of $E$, 
and let $\hh\colon \hX \to \tE$ be the pull-back of $h$ along $p$. 
We get an induced mapping, $\hat{p}\colon \hX \to X$, which 
is the Galois cover associated to the normal subgroup 
$\ker(h_{\sharp}) \subset \pi_1(X)$, as in diagram \eqref{eq:cd1}.
\begin{equation}
\label{eq:cd1}
\xymatrix{
H \ar^{\hat{i}}[r] \ar@{=}[d] 
& \hX\ar^\hh[r] \ar^{\hat{p}}[d] 
& \tE\ar^{p}[d]\\ 
H \ar^{i}[r] & X \ar^h[r] & E\\
H \ar^{i^*}[r] \ar@{=}[u] 
& X^* \ar^{h^*}[r] \ar[u]& E^* \ar_{\iota}[u]\\
}
\end{equation}

Note that $\tE$ is either the complex affine line $\C$ 
(when $\genus (E)=1$), or the open unit disc $\D$ 
(when $\genus (E)>1$).  Moreover,  $\hh\colon \hX \to \tE$ 
is a proper complex analytic mapping, with smooth fiber 
$H$.  Assuming $h$ has only isolated critical points, 
we infer that $\hh$ has countably many isolated singularities.

\begin{lemma}
\label{lem:xn}
The space $\hX$ is the union of an increasing sequence 
of open subsets $X_n$, such that, for each $n\ge 1$, 
the set $X_n$ contains $H$, and the inclusion 
$i_n\colon H \to X_n $ is an $(r-1)$-homotopy equivalence. 
\end{lemma}

\begin{proof}
Let $S=\tE \setminus \hh(C(\hh))$ be the set of regular values 
for $\hh$.  By applying a suitable automorphism to $\tE$, 
we may assume that $b=0$ belongs to $S$.
For a fixed $n$, define 
\[
X_n=\hh^{-1}(D_n),
\]
with $D_n$ an open disc in $\tE$ centered at $b$ and 
of radius $r_n=n$ (when $\tE=\C$), or $r_n=1-1/n$ 
(when $\tE=\D$).  Clearly, $X_n$ contains $H= \hh^{-1}(b)$.  

Consider a finite family of embeddings, 
$\gamma_c \colon [0,1] \to D_n$, parametrized 
by the critical values $c \in \hh (C(\hh))\cap D_n$, 
such that
\begin{romenum}
\item 
$\gamma_c(0)=b$, $\gamma_c(1)=c$, $\gamma_c((0,1)) \subset S$;
\item  
 $\gamma_c([0,1]) \cap \gamma_{c'}([0,1])=b$, for $c \ne c'$;
\item  
for each $c \in \hh (C(\hh))\cap D_n$, one can find 
a small closed disc $D_c\subset D_n$  centered at $c$, disjoint 
from the paths $\gamma_{c'}([0,1])$ for $c \ne c'$; 
\item 
$\gamma_c([0,1]) \cap \partial D_c=\gamma_c(1-\delta)$, for 
the same $\delta$, with $1\gg \delta >0$.
\end{romenum}

Let $K_n=\bigcup _c \gamma_c([0,1]) \cup \bigcup_c D_c$. 
Since $\hh$ is a fibration over $S$, it follows that
$X_n $ has the same homotopy type as $\hh^{-1}(K_n)$.
Similarly, let $L_n=\bigcup _c \gamma_c([0,1-\delta])$. 
Then $L_n$ is a contractible space and hence
$\hh^{-1}(L_n)$ has the homotopy type of $H= \hh^{-1}(b) $. 
Note that $\hh^{-1}(K_n)$ is obtained from $\hh^{-1}(L_n)$
by replacing the fibers $\hh^{-1}( \gamma_c(1-\delta)   )$ by the 
corresponding tubes $\hh^{-1}( D_c)$. Since we are in a proper 
situation, with finitely many isolated singularities, each such tube 
has the homotopy type of the central singular fiber $\hh^{-1}(c )$, 
which in turn is obtained from the nearby smooth fiber   
$\hh^{-1}( \gamma_c(1-\delta) )$ by attaching a finite number 
of $r$-dimensional cells.  (These cells are the cones over the 
corresponding $(r-1)$-dimensional vanishing cycles; see for 
instance the very similar proof in \cite[pp.~72--73]{L}).

The above argument shows that each $X_n$ has the homotopy 
type of a space obtained from the smooth fiber $H$ by attaching 
finitely many $r$-cells. Therefore, $\pi_{i}(X_n, H)=0$, 
for all $i<r$, and the conclusion follows.
\end{proof}

\subsection{}
\label{subs:proof thm C}
We are now ready to finish the proof of 
Theorem \ref{thm:cplxmorse}. 

Part \eqref{c1}.  From Lemma \ref{lem:xn}, and the fact 
that homotopy groups commute with direct limits, it 
follows that $\pi_{i}(\hX, H)=0$, for all $i<r$. 

Part \eqref{c2}. Since, by assumption, $r\ge 3$, the 
exact homotopy sequence of the pair $(\hX, H)$ 
shows that the inclusion $\hat{i}\colon H \to \hX $ 
induces an isomorphism on fundamental groups. 
Referring to diagram \eqref{eq:cd1}, it follows that 
$\ker(h_{\sharp})\cong \pi_1(\hX) \cong \pi_1(H)$.  
Combining this isomorphism with \eqref{eq:hsharp} 
yields the desired conclusion.  \hfill\qed

\begin{remark} 
\label{r1}
In the case when all the isolated singularities of $h$ are non-degenerate, 
the above proof essentially goes back to Lefschetz \cite{Lf}, see 
Lamotke \cite{La}, Section 5 and Section 8, particularly claim (8.3.2) 
and its proof.  The situation considered there corresponds to rational 
pencils, for which one may decompose the base of the pencil, 
$E=\PP^1$, as the union of two discs glued along their common 
boundary.  Hence, there is no need to pass to the universal 
cover $\tE$, as for irrational pencils, thus avoiding the difficulty 
of having to handle infinitely many critical values. 
\end{remark}

\section{Branched covers and elliptic pencils}
\label{sect:elliptic}

In this section, we construct a family of irrational pencils that 
satisfy the hypothesis of Theorem \ref{thm:cplxmorse}. To obtain 
these examples, we will replace the products of free groups 
on at least $2$ generators appearing within the framework of 
Bestvina--Brady theory (as in \S\ref{ss13}) by products of 
fundamental groups of smooth projective curves of genus 
at least $2$ (as in \S\ref{subs:prodcurves}).  

\subsection{}
\label{subsect:branched covers}
The starting point is a classical branched covering construction.
Let $E$ be an arbitrary complex elliptic curve. Let $B\subset E$ 
be a finite subset, of cardinality $\abs{B} =2g-2$, with $g>1$. Then
\begin{equation}
\label{eq:hsplit}
H_1(E\setminus B, \Z)\cong H_1(E, \Z) \oplus H_1^{B}\, ,
\end{equation}
where the group $H_1^{B}$ is generated by the homology classes 
$\{ \alpha_b \}_{b\in B}$ of elementary small positive loops around 
the points of $B$, subject to the single relation 
$\sum_{b\in B} \alpha_b =0$. 

Let $\varphi \in \Hom (\pi_1(E\setminus B), \Z/2 \Z)$ be any 
homomorphism with the property that, with respect to 
decomposition \eqref{eq:hsplit}, 
\begin{equation}
\label{eq:defphi}
\varphi (\alpha_b)=1,\ \forall b\in B. 
\end{equation}

The next result is of a well-known type. We include a 
sketch of proof, for the reader's convenience.

\begin{prop} 
\label{p1}
For any choice of $B\subset E$ and $\varphi$ as above, there 
is a projective, smooth curve $C$ of genus $g$, together with 
a ramified Galois $\Z/2\Z$-cover, $f\colon C \to E$. Furthermore, 
the map $f$ induces a bijection between the ramification locus 
$R\subset C$ and the branch locus $B\subset E$; the restriction 
$f\colon C\setminus R\to E\setminus B$ is the Galois cover 
corresponding to $\varphi$; and $f$ has ramification index $2$ 
at each point of $R$.
\end{prop}

\begin{proof}
Set $E^*=E\setminus B$, and let $f^*\colon C^*\to E^*$ be 
the Galois $\Z/2\Z$-cover associated to $\ker (\varphi)$. 
In view of a classical result of Stein \cite{St}, this cover 
extends uniquely to a ramified covering between the respective 
compactifications, $f\colon C \to E$. 

By construction, the ramification locus $R$ coincides 
with the critical set $C(f)$.  The assertions on the 
restriction of $f$ to $R$, and on ramification indices, 
are straightforward consequences of covering space theory. 
It follows that the ramification divisor of $f$ is $\sum_{c\in R} c$, 
whence $\genus (C)=g$, by the Riemann-Hurwitz formula; 
see for instance \cite[p.~142]{Mu}. 
\end{proof}

\subsection{}  
\label{subsec:elliptic curves}

Let $E$ be an elliptic curve, and fix an integer $r>1$. 
For each index $j$ from $1$ to $r$, let $B_j\subset E$ 
be a finite set with $\abs{B_j}= 2g_j -2>0$. Choose a 
homomorphism $\varphi_j \colon H_1(E\setminus B_j) 
\to \Z/2 \Z$ satisfying \eqref{eq:defphi}, and denote by 
$f_j\colon C_j\to E$ the corresponding branched cover 
(with ramification locus $R_j$ and branch locus $B_j$), 
as constructed in Proposition \ref{p1}.  Note that 
$\genus(C_j)=g_j\ge 2$. 

Write $X=C_1 \times \cdots \times C_r$  and 
$E^{\times r}=E\times \cdots \times E$, and 
consider the product mapping 
\[
f=f_1\times \cdots \times f_r\colon X \to E^{\times r}.
\]
Set $X_1= \prod_{j=1}^{r}( C_j \setminus R_j )$ and 
$Y_1=\prod_{j=1}^{r}(E \setminus B_j )$.
It follows that $f$ restricts to a $ (\Z/2\Z)^r$-covering, 
$f\colon X_1 \to Y_1$, which is determined by the 
homomorphism 
\begin{equation}
\label{eq:rcover}
\varphi:= \bigtimes_{j=1}^r \varphi_j \colon H_1(Y_1) \to (\Z/2\Z)^r .
\end{equation}

\subsection{}  
\label{subsec:elliptic pencil}
Let  $s_2\colon E^{\times 2} \to E$ be the group law of 
the elliptic curve, and extend it by associativity to a map 
$s_r\colon E^{\times r} \to E$.  Using these maps, we 
may define a holomorphic map $h$ as the composite 
\begin{equation}
\label{eq:hmap}
h=s_r \circ f\colon X \to E.
\end{equation}
The next result shows that $h$ is an {\em elliptic pencil}, that is,
an irrational pencil over an elliptic curve.

\begin{lemma} 
\label{lem3}
The smooth fiber of $h$ is connected.
\end{lemma}

\begin{proof}
Let $H_t=h^{-1}(t)$ be a smooth fiber of $h$.  
In order to show that $H_t$ is connected, 
it is enough to check that $H_1:=H_t \cap X_1$ is connected. 
This is due to the fact that no component of $H_t$ is contained 
in a hypersurface of the form $\{ (x_1,\dots, x_r)\in X \mid x_j=c \}$. 
Indeed, such an inclusion would force equality, whence 
$\sum_{i\ne j} f_i(x_i)= t- f_j(c)$, for all $x_i\in C_i$
($i\ne j$). Clearly, this is impossible, since $r\ge 2$.

Set $Z_t=s_r^{-1}(t)$. Note that $f\colon H_1 \to Z_t \cap Y_1$ 
is the pull-back of the covering $f\colon X_1 \to Y_1$, along the 
inclusion $\iota \colon Z_t \cap Y_1 \hookrightarrow Y_1$. 
Therefore, $H_1$ is connected, provided the composition 
$ H_1(Z_t \cap Y_1 ) \stackrel{\iota_*}{\longrightarrow} H_1(Y_1) 
\stackrel{\varphi}{\longrightarrow} (\Z/2\Z)^r$ is onto. To check this 
condition, it is enough to verify that each generator 
$\varepsilon_j \in (\Z/2\Z)^r$, with $1$ in the $j$-th position 
and $0$ elsewhere, lies in the image of $\varphi \circ \iota_*$.

Pick a point $b_j^0\in B_j$. We may then write 
the generic point $t\in E$ in the form 
$t=\sum_{i=1}^r t_i^0$, with $t_j^0= b_j^0$ and 
$t_i^0\notin B_i$, if $i\ne j$. 
Choose $i\ne j$ and define 
$\gamma =(\gamma_1,\dots, \gamma_r)\colon S^1\to E^r$
as follows: $\gamma_j$ is an elementary small positive loop 
around $b_j^0$;  $\gamma_i= t_i^0+ t_j^0- \gamma_j$; 
and $\gamma_k$ is the constant loop at $t_k^0$, for $k\ne i, j$. 
By our choices, $[\gamma]$ will be an element of 
$ H_1(Z_t \cap Y_1 )$, if $\gamma_j$ is small enough. 
Using \eqref{eq:defphi} and \eqref{eq:rcover}, it is readily 
seen that $\varphi \iota_*([\gamma])=\varepsilon_j$.
\end{proof}

\begin{lemma} 
\label{lem2}
The map $h\colon X \to E$ has only isolated singularities;  
more precisely, $C(h)=R_1\times \cdots\times R_r$.
Moreover, for each $p \in C(h)$, the induced function germ, 
$h\colon (X,p)\cong (\C^r,0) \to (E,h(p)) \cong (\C,0)$,
is a non-degenerate quadratic singularity, i.e., 
an $A_1$-singularity.  
\end{lemma}

\begin{proof}
Let $p=(p_1,\dots ,p_r) \in X=C_1 \times \cdots \times C_r$. 
Then, clearly, $\im d_pf=V_1 \times \cdots\times V_r$, 
where $V_j=0$ if $p_j \in C(f_j)$ and $V_j=\C$ otherwise. 
This implies the first claim.

The second claim follows from the fact that each function germ 
$f_j\colon(C_j,p_j)\cong (\C,0) \to (E,f_j(p_j))\cong (\C,0)$, 
where $p_j\in R_j$, is given in suitable coordinates 
by $x_j \mapsto x_j^2$, while the germ 
$s_r\colon (E^r,(f_1(p_1),\dots ,f_r(p_r)))\cong (\C^r,0) \to 
(E,h(p))\cong (\C,0)$ is given, again in suitable coordinates, 
by $(y_1,\dots ,y_r) \mapsto y_1+\cdots +y_r$.
\end{proof}

\begin{remark}
\label{rem:real Morse}
Denote by $H_t$ the generic smooth fiber of $h$, and 
assume $r\ge 3$.  Then $\pi_1(H_t) \cong \ker(h_{\sharp})$, 
as a consequence of Theorem 1.1 from Shimada \cite{Sh}, 
combined with our Lemmas \ref{lem3} and \ref{lem2}.
Of course, this also follows from Theorem \ref{thm:cplxmorse}, 
Part~\eqref{c2}.  

It is worth pointing out that our approach 
provides finer information, at the level of cell structures.
Indeed, let $h\colon X\to E$ be an elliptic pencil with only 
non-degenerate singularities, for instance, one of the pencils 
constructed above.  In this case, the function $g\colon \hX \to \RR$, 
given by $g(z)=\vert \hat{h} (z) \vert^2$, has only Morse singularities 
of index $r=\dim X$ on $\hX \setminus H_0$, as can be seen from 
the expansion 
\[
\Big( 1+\sum_{j=1}^{r} (x_j+ \sqrt{-1}\, y_j)^2 \Big) \cdot 
\Big( 1+\sum_{j=1}^{r} (x_j-\sqrt{-1}\, y_j)^2 \Big) = 
 1 + 2\sum_{j=1}^{r}(x_j^2-y_j^2)+ \cdots .
\]

Since $g$ is proper and the closed tube $T_0=g^{-1}([0,\epsilon])$ is 
homotopy equivalent to the central fiber $H_0=g^{-1}(0)$ for $\epsilon >0$ 
small enough, it follows from standard Morse theory (see \cite[I.3]{Mi}) 
that $\hX$ has the homotopy type of the smooth fiber $H_0$, with 
countably many $r$-cells attached. Clearly, these cells are indexed 
by $C(\hat{h})$, the set of critical points of $\hat{h}$.
\end{remark}

\section{Projective groups with exotic finiteness properties}
\label{sec:exotic}

In this section, we put things together, and finish the proof 
of Theorem \ref{thm:exko}.  

\subsection{}
\label{subs:bigthm}

We start by proving the following theorem. 

\begin{theorem}
\label{thm:exo}
Let $X$ be an irreducible, smooth projective variety of 
dimension $r \ge 3$.  Assume that  the universal cover 
of $X$ is an $(r-2)$-connected Stein manifold, 
and that $\V^r_1(\pi_1(X))= \T_{\pi_1(X)}$. 
Let $h\colon X\to E$ be a holomorphic map to an elliptic 
curve $E$, with connected smooth fiber $H$, and with isolated 
singularities. Then:

\begin{enumerate}
\item \label{e1}
The homotopy groups $\pi_i(H)$ vanish for $2\le i\le r-2$, 
while $\pi_{r-1}(H)\ne 0$.

\item \label{e0}
The universal cover $\tH$ is a Stein manifold.

\item \label{e2}
The group $N=\pi_1(H) $ has a $K(N, 1)$ with finite 
$(r-1)$-skeleton, but $H_r(N, \Z)$ is not finitely generated.

\item \label{e3}
The group $N$ is not commensurable (up to finite kernels) 
to any group $N'$ having a $K(N', 1)$ of finite type.
\end{enumerate}
\end{theorem}

\begin{proof}
As before, let $p\colon\tE \to E$ be the universal cover, and 
let $\hh\colon \hX\to \tE$ be the pull-back of $h$ along $p$. 
The universal cover $\tX \to X$ factors as $\hat{p}\circ q$, 
where $\hat{p}\colon \hX\to X$ is the pull-back of $p$ along $h$, 
and $q\colon \tX \to \hX$ is the universal cover of $\hX$. 

Part \eqref{e1}. The vanishing property for the higher homotopy 
groups is a consequence of Theorem \ref{thm:cplxmorse}, given 
our connectivity assumptions on the universal cover $\tX$. 
If $\pi_{r-1}(H)$ would also vanish, we could construct a 
classifying space $K(N, 1)$ by attaching to $H$ cells of 
dimension $r+1$ and higher. In particular, $N$ would be 
of type $\F_r$, with finitely generated $r$-th homology group, 
contradicting property \eqref{e2}, which is proved below.

Part \eqref{e0}. Since the pair $(\hX, H)$ is $(r-1)$-connected 
and $r\ge 3$, the universal cover $\tH$ coincides with $q^{-1}(H)$. 
This is a closed analytic submanifold of the Stein manifold $\tX$, 
hence Stein as well.

Part \eqref{e2}. The group  $N$ is of type $\F_{r-1}$, by the 
same argument as in the proof of Part \eqref{e1}. Assuming 
$H_r(N, \Z)$ to be finitely generated, we infer that 
$\dim_{\C} H_{\le r}(N, \C)<\infty$.

On the other hand, we know from Theorem \ref{thm:cplxmorse} 
that $h_{\sharp}\colon \pi_1(X) \to \pi_1(E)=\Z^2$ is surjective, 
and $N\cong \ker (h_{\sharp})$. But this contradicts 
Theorem \ref{thm:vfin}, due to our hypothesis on 
$\V^r_1(\pi_1(X) )$. 

Part \eqref{e3}. Follows from Part \eqref{e2} and 
Proposition \ref{prop:commens}.
\end{proof}

\subsection{}
\label{subs:highpi}
As a by-product of our Morse-theoretical approach, we can 
give a precise description of both the $\Z N$-structure of 
$\pi_{r-1}(H)$ and the homotopy type of $\tH$, in \eqref{e1} 
and \eqref{e0} above. We keep the notation and hypothesis 
from Theorem \ref{thm:exo}. 

\begin{corollary}
\label{cor:hipi}
Assume that the universal cover of $X$ is $r$-connected, 
and the map $h\colon X\to E$ has only non-degenerate 
singularities. Let $C(h)$ be the set of critical points of $h$, 
and let $H$ be the smooth fiber of $h$.   Then:
\begin{enumerate}
\item \label{hi1}
The homotopy group  $\pi_{r-1}(H)$ is a free $\Z \pi_1(H)$-module, 
with generators in one-to-one correspondence with $C(h)\times \pi_1(E)$.

\item \label{hi2}
The universal cover $\tH$ has the homotopy type of a wedge of 
$(r-1)$-spheres, indexed by $\pi_1(H) \times C(h)\times \pi_1(E)$.
\end{enumerate}
\end{corollary}

\begin{proof}
Part \eqref{hi1}. The exact homotopy sequence of the pair 
$(\hX, H)$ identifies $\pi_{r-1}(H)$ with $\pi_r(\hX, H)$, as 
modules over $\Z \pi_1(H)$. Now recall from 
Remark \ref{rem:real Morse} that $\hX$ has the homotopy 
type of $H$, with some $r$-cells attached. Moreover, these 
cells are indexed by  
$C(\hh)= \hat{p}^{-1}(C(h))\cong C(h)\times \pi_1(E)$.
Since $r\ge 3$, the claim follows from \cite[Exercise 7.F.3]{Sp}.

Part \eqref{hi2}. We know from Theorem \ref{thm:exo}\eqref{e0} 
that $\tH$ is an $(r-1)$-dimensional Stein manifold. Therefore, 
$\tH$ has the homotopy type of a CW-complex of 
dimension at most $r-1$. 
Let $W$ be the bouquet of $(r-1)$-spheres indexed by 
$\pi_1(H) \times C(h)\times \pi_1(E)$, and let 
$\psi \colon W \to \tH$ be the map whose restriction 
to each sphere represents the corresponding generator 
of the free abelian group $\pi_{r-1}(\tH)$, computed 
in Part \eqref{hi1}.  By Theorem \ref{thm:exo}\eqref{e1}, 
$\pi_{i}(\tH)= \pi_{i}(W)=0$, for $1\le i\le r-2$.  Moreover, 
$\psi$ induces an isomorphism on $\pi_{r-1}$, by construction. 
The claim follows from the Hurewicz and Whitehead theorems.
\end{proof}

Similar highly connected Stein spaces having the homotopy type
of bouquets of spheres occur in the study of local complements
of isolated non-normal crossing singularities, see \cite{Li}, \cite{DLi}.

\subsection{}
\label{subs:finish thm A}
We can now finish the proof of Theorem \ref{thm:exko}. 

Fix an integer $r\ge 3$, and let $E$ be an elliptic curve. 
For each index $j$ from $1$ to $r$, construct a 
$2$-fold branched cover $f_j\colon C_j \to E$, 
with $C_j$ a curve of genus at least $2$,  
as in \S\ref{subsec:elliptic curves}. The product 
$X=\prod_{j=1}^r C_j$ is a smooth, projective variety 
of dimension $r$, whose universal cover is a contractible, 
Stein manifold.  By Corollary \ref{cor:product of curves}\eqref{cc1}, 
the characteristic variety $\V^r_1(\pi_1(X))$ coincides 
with the character torus $\T_{\pi_1(X)}$.  

Now define a holomorphic map $h\colon X\to E$ as 
in \eqref{eq:hmap}.  By Lemmas \ref{lem3} and \ref{lem2}, 
the smooth fiber of $h$ is connected, and $h$ has only 
isolated, non-degenerate singularities.   Thus, the 
hypotheses of Theorem \ref{thm:exo} hold.  The 
conclusions of Theorem \ref{thm:exko} follow at 
once from Theorem \ref{thm:exo}. \hfill\qed
\smallskip

For this class of examples, the conclusions of 
Corollary \ref{cor:hipi} are valid as well.

\subsection{}
\label{subs:cdim}

The Stein condition influences another finiteness property 
of projective groups, namely, their cohomological dimension. 

\begin{prop}
\label{prop:cdim}
Let $M$ be a compact connected, $m$-dimensional complex analytic 
manifold, and let $G=\pi_1(M)$. If the universal cover $\tM$ 
is Stein, then $\cd (G)\ge m$. 
\end{prop}

\begin{proof}
Let $\kappa \colon M\to K(G, 1)$ be a classifying map, with homotopy 
fiber $\tM$. Let us examine the associated Serre spectral sequence,
\[
E^2_{st}= H_s(G, H_t(\tM, \Z)) \Rightarrow H_{s+t}(M, \Z)\, .
\]
Since $\tM$ is Stein, it has the homotopy type of a CW-complex 
of dimension at most $m$. Therefore, $E^2_{st}=0$, for $t>m$.
Assuming $\cd (G)<m$, we infer that $E^2_{st}=0$, for $s\ge m$. 
These two facts together imply that $H_{2m}(M, \Z)=0$, a contradiction.
\end{proof}

A related statement holds for the smooth fiber $H$ from 
Theorem \ref{thm:exo}: $\cd (\pi_1(H))\ge \dim H +1$, 
as follows from Part \eqref{e2}.

\begin{ack}
We are grateful to J\'{a}nos Koll\'{a}r for bringing up to our attention 
reference \cite{K}, and for stimulating our interest in finding a projective 
analog of our results on Bestvina-Brady groups from \cite{DPS06}.
\end{ack}

\bibliographystyle{amsplain}

\end{document}